# Existence and uniqueness for a class of nonlinear higher-order partial differential equations in the complex plane

O. Costin[*] and S. Tanveer[†]

November 5, 2018

### Abstract

We prove existence and uniqueness results for nonlinear third order partial differential equations of the form

$$f_t - f_{yyy} = \sum_{j=0}^{3} b_j(y,t;f)\ f^{(j)} + r(y,t)$$

where superscript $j$ denotes the $j$-th partial derivative with respect to $y$. The inhomogeneous term $r$, the coefficients $b_j$ and the initial condition $f(y,0)$ are required to vanish algebraically for large $|y|$ in a wide enough sector in the complex $y$-plane. Using methods related to Borel summation, a unique solution is shown to exist that is analytic in $y$ for all large $|y|$ in a sector. Three partial differential equations arising in the context of Hele-Shaw fingering and dendritic crystal growth are shown to be of this form after appropriate transformation, and then precise results are obtained for them. The implications of the rigorous analysis on some similarity solutions, formerly hypothesized in two of these examples, are examined.

## 1 Introduction

The theory of partial differential equations (PDEs), when one or more of the independent variables are in the complex plane, appears to be largely undeveloped. The classic Cauchy-Kowaleski theorem holds for a system of first-order equations (or those equivalent to it) when the quasi-linear equations have analytic coefficients and analytic initial data is specified on an analytic but non-characteristic curve. Then, the C-K theorem guarantees the local existence and uniqueness of analytic solutions. As is well known, its proof relies on the convergence of local power series expansions and, without the given hypotheses, the power series may have zero radius of convergence and the C-K method does not yield solutions.

---
[*]Mathematics Department, Rutgers University, Busch Campus, Hill Center, 110 Frelinghuysen Rd., Piscataway NJ 08854; costin@math.rutgers.edu; http://www.math.rutgers.edu/∼costin.
[†]Mathematics Department, The Ohio State University, 231 W. 18th Avenue, Columbus, OH 43210; tanveer@math.ohio-state.edu; http://www.math.ohio-state.edu/∼tanveer.

Relatedly, not much is known in general for higher-order nonlinear partial differential equations in the complex plane, even for analytic initial conditions and analytic dependence of the coefficients. The only work we are aware of on nonlinear partial differential equations in the complex plane involving higher spatial derivatives is that of Sammartino & Caflisch [1, 2] who proved among other results the existence of a solution to nonlinear Prandtl boundary layer equations for analytic initial data. This work involved inversion of the heat operator $\partial_t - \partial_{YY}$ and using the abstract Cauchy-Kowalewski theorem for the resulting integral equation. Unfortunately, this methodology cannot be adapted to our problem. The coefficients of the highest (third order) spatial derivatives in our equation depend on the unknown function as well. These terms cannot be controlled by inversion of a linear operator and estimates of the kernel, as used by Sammartino & Caflisch. Instead, the essence of the methodology introduced here is the use of large $y$ asymptotics, conveniently expressed in terms of the behavior of the unknown function in the Borel transform variable $p$ for small $p$. The choice of appropriate Banach spaces proves to be crucial, and after this choice the contraction mapping argument itself is not difficult.

One aim of the present paper is to obtain actual solutions with good smoothness and asymptotic properties for a class of PDEs, when power series solutions may have zero radius of convergence. Our approach based on Borel summation techniques provides at the same time appropriate existence and uniqueness results for a class of nonlinear PDEs in the complex domain.

Keeping in mind applications, we develop the framework for certain higher order partial differential equations in a domain where one of the independent variables ($y$ in this case) is complex, while the other ($t$) is real. While more sweeping generalizations are under way, the current paper is restricted in scope by the applications we have in mind and simplicity of exposition.

There is a class of nonlinear PDEs that have recently arisen in applications. The basic features of the application problems is that in the absence of a regularization (like surface tension) the initial value problem in the real domain is relatively simple; yet ill-posed in the sense of Hadamard for any Sobolev norm on the real domain. However, the analytically continued equations into the complex spatial domain are well-posed, even without a regularization term. Earlier, Garabedian [4] had recognized the conversion of ill-posed elliptic initial value problem into a well-posed one by excursion into the complex plane in the spatial variable. Moore [5], [6], Caflisch & Orellana [7], [8], Caflisch & Semmes [9] and Caflisch *et al* [10] have studied solutions to the complex plane equations that arose from simplifications of vortex sheet evolution (in fluid mechanical contexts). The initial value problem in these cases is ill-posed in the real domain, though well-posed in an appropriate class of analytic functions in a domain in the complex plane. Study of the complex equations proved useful since evidence [11] of finite time singularities in the real domain can be traced to earlier singularity formation in the complex domain.

In the physical context of Hele-Shaw dynamics, it was suggested [12] that it is fruitful to study the complex plane equations even when the initial value problem is well-posed in the physical domain through the addition of a small regularization term. The advantage of this procedure is that one can study small regularization effects by perturbing about the relatively simpler but well-posed zeroth order problem. The ill-posedness of the unregularized problem in the real domain, shown earlier by Howison [15], is transferred into ill-posedness of the analytic continuation of initial data to the complex plane. However, when analytic initial data is specified in a domain in the complex plane such as to allow for isolated singularities, there is no ill-posedness of the zeroth order approximation of the dynamics. This provides the basis for a perturbative study that includes small but nonzero regularization effects in the real domain. Consideration of an ensemble of complex initial conditions, subject to appropriate constraints on its behavior on the real axis, provide a way to understand the robust features of the dynamics when regularization effects are small. Indeed this procedure has yielded information about how small surface tension can singularly perturb a smooth solution of the un-regularized dynamics [13]. It has given scaling results on nonlinear dendritic processes as well [17]. However, much of the results derived so far are purely formal and rely fundamentally on the existence and uniqueness of analytic solutions to certain higher order nonlinear partial differential equations in a sector in the complex plane, with imposed far-field matching conditions. Indeed, in a more general context, one can expect that whenever regularization appears in the form of a small coefficient multiplying the highest spatial derivative, the resulting asymptotic equation in the neighborhood of initial complex singularities will satisfy higher order nonlinear partial differential equation with sectorial far-field matching condition in the complex plane of the type shown in Examples 1-3. Hence, there is a need to develop a general theory in this direction.

## 2    Problem Statement and Main Result

We seek to prove the existence and uniqueness of solutions $f(y,t)$ to the initial value problem for a general class of quasilinear partial differential equations of the form:

$$f_t - f_{yyy} = \sum_{j=0}^{3} b_j(y,t;f)\, f^{(j)} + r(y,t) \quad \text{with} \quad f(y,0) = f_I(y) \tag{1}$$

where the superscript $^{(j)}$ refers to the $j$-th derivative with respect to $y$. The inhomogeneous term $r(y,t)$ is a specified analytic function in the domain

$$\mathcal{D}_{\rho_0} = \{(y,t) : \arg y \in (-2\pi/3, 2\pi/3),\ |y| > \rho_0 > 0, 0 \le t \le T\}$$

and it is assumed that in $\mathcal{D}_{\rho_0}$ there exist constants $\alpha_r \ge 1$ and $A_r$, with only $A_r$ allowed to depend on $T$, such that

$$|y^{\alpha_r}\, r(y,t)| < A_r(T) \tag{2}$$

Further in (1), the coefficients $b_j$ may depend on the solution $f$– this is how nonlinearity in the problem arises. It is possible to extend the current theory to include dependence of $b_j$ on $f^{(j)}$ as well, though for simplicity we will restrict only to dependence on $f$. Further, we restrict to the case where each $b_j$ is given by a convergent series

$$b_j(y,t;f) = \sum_{k=0}^{\infty} b_{j,k}(y,t)\, f^k \tag{3}$$

for known $b_{j,k}$, analytic for $y$ in $\mathcal{D}_{\rho_0}$. It will be assumed that in this domain, there exists some choice of positive constants $\beta$, $\alpha_j$, and $A_b$, independent of $j$ and $k$ (with $\beta$ and $\alpha_j$ independent of $T$ as well), such that

$$|y^{\alpha_j + k\beta}\, b_{j,k}| < A_b(T) \tag{4}$$

Further, the series (3) converges in the domain $\mathcal{D}_{\phi,\rho}$, defined as

$$\mathcal{D}_{\phi,\rho} = \left\{(y,t) : \arg y \in \left(-\frac{\pi}{2} - \phi, \frac{\pi}{2} + \phi\right), \ |y| > \rho > \rho_0 \right.$$
$$\left. \text{where } 0 < \phi < \frac{\pi}{6}, \ 0 \leq t \leq T \right\} \quad (5)$$

if

$$|f| < \rho^\beta \quad (6)$$

**Condition 1** *The solution $f(y,t)$ we seek for (1) is required to be analytic for complex $y$ in $\mathcal{D}_{\phi,\rho}$ for some $\rho > 0$ (to be determined later). In the same domain, the solution and the initial condition $f_I(y)$ must satisfy the condition*

$$|y^{\alpha_r} f(y,t)| < A_f(T) \quad (7)$$

*for some $A_f$ that can only depend on $T$ for $(y,t) \in \mathcal{D}_{\phi,\rho}$.*

It is clear that for large $y$ such a solution $f$ will indeed satisfy (6), the condition for the convergence of the infinite series in (1). The general theorem proved in this paper is that

**Theorem 2** *For any $T > 0$ and $0 < \phi < \pi/6$, there exists $\tilde{\rho}$ such that the partial differential equation (1) has a unique solution $f$ that is analytic in $y$ and $O(y^{-1})$ as $y \to \infty$ for $(y,t) \in \mathcal{D}_{\tilde{\rho},\phi}$. In fact, we have for this solution $f = O(y^{-\alpha_r})$ as $y \to \infty$.*

We note that uniqueness requires analyticity and decay properties of $y$ in a large enough sector. The proof of theorem (2) will have to await some definitions and lemmas. It is to be noted that, from a formal argument if $f$ is small, the dominant balance for large $y$ is between $f_t$ on the left of (1) and $r(y,t)$ on the right, indicating that $f(y,t) \sim f_I(y) + \int_0^t r(y,t) \, dt$. Since each of $f_I(y)$ and $r(y,t)$ decays algebraically as $y \to \infty$ within $\mathcal{D}_{\rho_0}$ at a rate $y^{-\alpha_r}$, which for $\alpha_r > 1$ is much less than $y^{-1}$, this suggests suggests that other terms in the differential equation (1) should not contribute. This is in fact shown rigorously for $|y|$ large with $\arg y \in (-2\pi/3, 2\pi/3)$. As shall be seen in the examples, this behavior for solution $f$ is not valid outside this sector, where in general, one can expect infinitely many singularities with an accumulation point at $\infty$.

## 3 Inverse Laplace Transform and Equivalent Integral Equation

The Inverse Laplace Transform $G(p,t)$ of a function $g(y,t)$ analytic in $\mathcal{D}_{\phi,\rho}$ and vanishing algebraically as $|y| \to \infty$ is given by:

$$G(p,t) = \left[\mathcal{L}^{-1}\{g\}\right](p,t) \equiv \frac{1}{2\pi i} \int_{\mathcal{C}_D} e^{py} g(y,t) \, dy \quad (8)$$

where $\mathcal{C}_D$ is a contour as in Fig. 1 (or deformations thereof), entirely within the domain $\mathcal{D}_{\phi,\rho}$. We restrict $p$ to the domain:

$$\mathcal{S}_\phi \equiv \{p : \arg p \in (-\phi, \phi), 0 < |p| < \infty \}$$

It is easily seen that if $g(y,t) = y^{-\alpha}$ for $\alpha > 0$, then $G(p,t) = p^{\alpha-1}/\Gamma(\alpha)$. From the following Lemma, it is clear that the same kind of behavior for the Inverse Laplace Transform $G(p,t)$ can be expected for small $p$ in $\mathcal{S}_\phi$, with a $y^{-\alpha}$ behavior of $g$ at $\infty$.

**Lemma 3** *If $g(y,t)$ is analytic in $y$ in $\mathcal{D}_{\phi,\rho}$, and satisfies*

$$|y^\alpha \, g(y,t)| < A(T) \tag{9}$$

*for $\alpha \geq \alpha_0 > 0$, then for any $\delta \in (0, \phi)$ the Inverse Laplace Transform $G = \mathcal{L}^{-1} g$ exists in $\mathcal{S}_{\phi-\delta}$ and satisfies*

$$|G(p,t)| < C \frac{A(T)}{\Gamma(\alpha)} |p|^{\alpha-1} e^{2|p|\rho} \tag{10}$$

*for some $C = C(\delta, \alpha_0)$.*

PROOF. We first consider the case when $2 \geq \alpha \geq \alpha_0$. Let $C_{\rho_1}$ be the contour $C_D$ in Fig. 1 that passes through the point $\rho_1 + |p|^{-1}$, and given by $s = \rho_1 + |p|^{-1} + ir \exp(i\phi \operatorname{signum}(r))$ with $r \in (-\infty, \infty)$. Choosing $2\rho > \rho_1 > (2/\sqrt{3})\rho$, we have $|s| > \rho$ along the contour and therefore, with $\arg(p) = \theta \in (-\phi + \delta, \phi - \delta)$,

$$|g(s,t)| < A(T)|s|^{-\alpha} \quad \text{and} \quad |e^{sp}| \leq e^{\rho_1 |p|+1} e^{-|r||p|\sin|\phi-\theta|}$$

Thus

$$\left| \int_{C_{\rho_1}} e^{sp} g(s,t) ds \right| \leq 2A(T) e^{\rho_1 |p|+1} \int_0^\infty \left| \rho_1 + |p|^{-1} + ire^{i\phi} \right|^{-\alpha} e^{-|p|r \sin \delta} dr$$

$$\leq \tilde{K} A(T) e^{\rho_1 |p|} \left| \rho_1 + |p|^{-1} \right|^{-\alpha} \int_0^\infty e^{-|p|r \sin \delta} dr \leq K \delta^{-1} |p|^{\alpha-1} e^{2\rho|p|} \tag{11}$$

where $\tilde{K}$ and $K$ are constants independent of any parameter. Thus, the Lemma follows for $2 \geq \alpha \geq \alpha_0$, if we note that $\Gamma(\alpha)$ is bounded in this range of $\alpha$, with the bound only depending on $\alpha_0$.

For $\alpha > 2$, there exists an integer $k > 0$ so that $\alpha - k \in (1,2]$. Taking $(k-1)! h(y,t) = \int_\infty^y g(z,t)(y-z)^{k-1} dz$ (clearly $h$ is analytic in $\mathcal{D}_{\phi,\rho}$ and $h^{(k)}(y,t) = g(y,t)$), we get

$$h(y,t) = \frac{(-y)^k}{(k-1)!} \int_1^\infty g(yp,t)(p-1)^{k-1} dp$$

$$= \frac{(-1)^k y^{k-\alpha}}{(k-1)!} \int_1^\infty A(yp,t) p^{-\alpha}(p-1)^{k-1} dp$$

with $|A(yp,t)| < A(T)$, whence

$$|h(y,t)| < \frac{A(T)\Gamma(\alpha-k)}{|y|^{\alpha-k}\Gamma(\alpha)}$$

From what has been already proved, with $\alpha - k$ playing the role of $\alpha$,

$$|\mathcal{L}^{-1}\{h\}(p,t)| < C(\delta)\frac{A(T)}{\Gamma(\alpha)}|p|^{\alpha-k-1}e^{2|p|\rho}$$

Since $G(p,t) = (-1)^k p^k \mathcal{L}^{-1}\{h\}(p,t)$, by multiplying the above equation by $|p|^k$, the Lemma follows for $\alpha > 2$ as well.
□

**Comment 1:** As mentioned before, when $g(y,t) = A(t)y^{-\alpha}$ we have $G(p,t) = \frac{A(t)}{\Gamma(\alpha)}p^{\alpha-1}$. In this case, the exponential factor in (10) can be omitted because of the algebraic behavior of $g(y,t)$ for all $y$. This result is relevant for Examples 1-3.

**Comment 2:** The constant $2\rho$ in the exponential bound can be lowered to anything exceeding $\rho$, but (10) suffices for our purposes.

**Comment 3:** Corollary 4 below implies that for any $p \in \mathcal{S}_\phi$, the Inverse Laplace Transform exists for the specified functions $r(y,t)$, $b_{j,k}(y,t)$, as well as the solution $f(y,t)$ to (1), whose existence is shown in the sequel.

**Comment 4:** Conversely, if $G(p,t)$ is any integrable function satisfying the exponential bound in (10), it is clear that the Laplace Transform along a ray

$$\mathcal{L}_\theta G \equiv \int_0^{\infty e^{i\theta}} dp \; e^{-py} \; G(p,t) \tag{12}$$

exists and defines an analytic function of $y$ in the half-plane $\Re[e^{i\theta}y] > 2\rho$ for $\theta \in (-\phi, \phi)$.

**Comment 5:** The next corollary shows that there exist bounds for $\mathcal{L}^{-1}\{b_{j,k}\}$ and $\mathcal{L}^{-1}\{r\}$ independent of arg $p$ in $\mathcal{S}_\phi$, because of the assumed analyticity and decay properties in the region $\mathcal{D}_{\rho_0}$, which contains $\mathcal{D}_{\phi,\rho}$.

**Corollary 4** *The Inverse Laplace Transform of the coefficient functions $b_{j,k}$ and the inhomogeneous function $r(y,t)$ satisfy the following upper bounds for any $p \in \mathcal{S}_\phi$*

$$|B_{j,k}(p,t)| < \frac{C_1(\phi,\alpha_j)}{\Gamma(\alpha_j+k\beta)} A_b(T) |p|^{k\beta+\alpha_j-1} e^{2\rho_0|p|} \tag{13}$$

$$|R(p,t)| < \frac{C_2(\phi)}{\Gamma(\alpha_r)} A_r(T) |p|^{\alpha_r-1} e^{2\rho_0|p|} \tag{14}$$

PROOF. From the assumed conditions we see that $b_{j,k}$ is analytic in $y$ over $\mathcal{D}_{\phi_1,\rho_0}$ for any $\phi_1$ satisfying $\pi/6 > \phi_1 > \phi$. So Lemma 3 can be applied for $g(y,t) = b_{j,k}$, with $\phi_1 = \phi + (\pi/6 - \phi)/2$ replacing $\phi$, and with $\delta$ replaced by $\phi_1 - \phi = (\pi/6 - \phi)/2$, and the same applies to $R(p,t)$, leading to (13) and (14). In the latter case, since $\alpha_r \geq 1$, $\alpha_0$ in Lemma 3 can be chosen to be 1. Thus, one can choose $C_2$ to be independent of $\alpha_r$, as indicated in (14). □

The formal Inverse Laplace Transform of (1) with respect to $y$ is

$$F_t + p^3 F = \sum_{j=0}^{3} (-1)^j \sum_{k=0}^{\infty} \left[ B_{j,k} * (p^j F) * F^{*k} \right](p,t) + R(p,t) \tag{15}$$

where the symbol $*$ stands for convolution (see also [3]). On formally integrating (15) with respect to $t$, we obtain the integral equation

$$F(p,t) = \sum_{j=0}^{3} \sum_{k=0}^{\infty} \int_0^t (-1)^j e^{-p^3(t-\tau)} \left[ (p^j F) * B_{j,k} * F^{*k} \right](p,\tau) d\tau + F_0(p,t)$$
$$\equiv \mathcal{N} F(p,t) \tag{16}$$

where

$$F_0(p,t) = e^{-p^3 t} F_I(p) + \int_0^t e^{-p^3(t-\tau)} R(p,\tau) \, d\tau \tag{17}$$

Here $F_I = \mathcal{L}^{-1}\{f_I\}$.

Our strategy is to reduce the problem of existence and uniqueness of a solution of (1) to the problem of existence and uniqueness of a solution of (16), under appropriate conditions.

## 4  Solution to the Integral Equation

To establish the existence and uniqueness of solutions to the integral equation, we need to introduce an appropriate function class for both the solution and the coefficient functions.

**Definition 5** *Denoting by $\overline{\mathcal{S}_\phi}$ the closure of $\mathcal{S}_\phi$, $\partial \mathcal{S}_\phi = \overline{\mathcal{S}_\phi} \setminus \mathcal{S}_\phi$ and $\mathcal{K} = \overline{\mathcal{S}_\phi} \times [0,T]$, we define for $\nu > 0$ (later to be taken appropriately large) the norm $\| \cdot \|_\nu$ as*

$$\|G\|_\nu = M_0 \sup_{(p,t) \in \mathcal{K}} (1+|p|^2)\, e^{-\nu |p|} \, |G(p,t)| \tag{18}$$

*where*

$$M_0 = \sup_{s \geq 0} \left\{ \frac{2(1+s^2)\left(\ln(1+s^2) + s \arctan s\right)}{s(s^2+4)} \right\} = 3.76 \cdots \tag{19}$$

**Definition 6** *We now define the following class of functions:*

$$\mathcal{A}_\phi = \{F : F(\cdot, t) \text{ analytic in } \mathcal{S}_\phi \text{ and continuous in }$$

$$\overline{\mathcal{S}_\phi} \text{ for } t \in [0, T] \text{s.t. } \|F\|_\nu < \infty\}$$

*It is clear that $\mathcal{A}_\phi$ forms a Banach space.*

**Comment 6:** Note that given $G \in \mathcal{A}_\phi$, $g(y,t) = \mathcal{L}_\theta\{G\}$ exists for appropriately chosen $\theta$ when $\rho$ is large enough so that $\rho \cos(\theta + \arg y) > \nu$, and that $|y\, g(y,t)|$ is bounded for $y$ on any fixed ray in $\mathcal{D}_{\phi,\rho}$.

**Lemma 7** *For $\nu > 4\rho_0 + \alpha_r$, $F_I = \mathcal{L}^{-1}\{f_I\}$ satisfies*

$$\|F_I\|_\nu < C(\phi) A_{f_I} (\nu/2)^{-\alpha_r+1}$$

*while $R = \mathcal{L}^{-1} r$ satisfies the relation*

$$\|R\|_\nu < C(\phi) A_r(T) (\nu/2)^{-\alpha_r+1}$$

*and therefore*

$$\|F_0\|_\nu < C(\phi)(T\, A_r + A_{f_I})(\nu/2)^{-\alpha_r+1} \tag{20}$$

PROOF. First note the bounds on $R$ in Corollary 4. We also note that $\alpha_r \geq 1$ and that for $\nu > 4\rho_0 + \alpha_r$ we have

$$\sup_p \frac{|p|^{\alpha_r-1}}{\Gamma(\alpha_r)} e^{-(\nu-2\rho_0)|p|} \leq \frac{(\alpha_r-1)^{\alpha_r-1}}{\Gamma(\alpha_r)} e^{-\alpha_r+1} (\nu-2\rho_0)^{-\alpha_r+1}$$

$$\leq K \alpha_r^{-1/2} (\nu/2)^{-\alpha_r+1}$$

where $K$ is independent of any parameter. The latter inequality follows by accounting for Stirling's formula for $\Gamma(\alpha_r)$ for large $\alpha_r$. Similarly,

$$\sup_p \frac{|p|^{\alpha_r+1}}{\Gamma(\alpha_r)} e^{-(\nu-2\rho_0)|p|} \leq \frac{(\alpha_r+1)^{\alpha_r+1}}{\Gamma(\alpha_r)} e^{-\alpha_r-1} (\nu-2\rho_0)^{-\alpha_r-1}$$

$$\leq K \alpha_r^{3/2} (\nu/2)^{-\alpha_r-1}$$

Using the definition of the $\nu-$norm and the two equations above, the inequality for $\|R\|_\nu$ follows. Since $f_I(y)$ is required to satisfy the same bounds as $r(y,t)$, a similar inequality holds for $\|F_I\|_\nu$. Now, from the relation (17),

$$|F_0(p,t)| < |F_I(p)| + T \sup_{0 \leq t \leq T} |R(p,t)|$$

Therefore, (20) follows. □

**Comment 7**: Not all Laplace-transformable analytic functions in $\mathcal{D}_{\phi,\rho}$ belong to $\mathcal{A}_\phi$. For the applications we have in mind, the coefficients are not bounded near $p = 0$ and hence do not belong in $\mathcal{A}_\phi$. It is then useful to introduce the following function class:

**Definition 8**

$$\mathcal{H} \equiv \left\{ H : H(p,t) \text{ analytic in } \mathcal{S}_\phi, |H(p,t)| < C|p|^{\alpha-1} e^{\rho|p|} \right\}$$

*for some positive constants $C$ and $\alpha$ and $\rho$ which may depend on $H$.*

**Lemma 9** *If $H \in \mathcal{H}$ and $F \in \mathcal{A}_\phi$, then for $\nu > \rho + 1$, $H * F$ belongs to $\mathcal{A}_\phi$, and satisfy the following inequality:*[1]

$$\|H * F\| \le \||H| * |F|\|_\nu \le C\, \Gamma(\alpha)(\nu - \rho)^{-\alpha} \, \|F\|_\nu \tag{21}$$

PROOF. From elementary properties of convolution, it is clear that $H * F$ is analytic in $\mathcal{S}_\phi$ and is continuous on $\overline{\mathcal{S}_\phi}$. Let $\theta = \arg\, p$. It is to be noted that

$$|H * F(p)| \le ||H| * |F||(p) \le \int_0^{|p|} |H(se^{i\theta})||F(p - se^{i\theta})|ds$$

But

$$|H(se^{i\theta})| \le Cs^{\alpha-1} e^{|s|\rho}$$

And

$$\int_0^{|p|} s^{\alpha-1} e^{|s|\rho} |F(p - se^{i\theta})| ds \le \|F\|_\nu e^{\nu|p|} |p|^\alpha \int_0^1 \frac{s^{\alpha-1} e^{-(\nu-\rho)|p|s}}{M_0(1 + |p|^2(1-s)^2)} ds \tag{22}$$

If $|p|$ is large, noting that $\nu - \rho \ge 1$, we obtain from Watson's lemma,

$$\int_0^{|p|} s^{\alpha-1} e^{|s|\rho} |F(p - se^{i\theta})| ds \le K\Gamma(\alpha)\|F\|_\nu \frac{e^{\nu|p|}}{M_0(1+|p|^2)} |\nu - \rho|^{-\alpha} \tag{23}$$

Now, for any other $|p|$, we obtain from (22),

$$\int_0^{|p|} s^{\alpha-1} e^{|s|\rho} |F(p - se^{i\theta})| ds \le K|\nu - \rho|^{-\alpha} \|F\|_\nu \frac{e^{\nu|p|}\Gamma(\alpha)}{M_0}$$

Thus the relation (23) must hold in general as it subsumes the above relation when $|p|$ is not large. From this relation, (21) follows from applying the definition of $\|.\|_\nu$.
□

**Corollary 10** *For $F \in \mathcal{A}_\phi$, and $\nu > 4\,\rho_0 + 1$,*

$$\|B_{j,k} * F\|_\nu \le \||B_{j,k}| * |F|\|_\nu \le K\, C_1(\phi, \alpha_j)\, (\nu/2)^{-k\beta - \alpha_j} A_b(T)\, \|F\|_\nu$$

PROOF. The proof follows simply by using Lemma 9, with $H$ replaced by $B_{j,k}$ and using the relations in Corollary 4. □

---

[1] In the following equation, $\| \cdot \|_\nu$ is extended naturally to continuous functions in $\mathcal{K}$.

**Lemma 11** *For $F, G \in \mathcal{A}_\phi$ and $j \geq 0$,*

$$|(p^j\ F) * G(p,t)| \ \leq\ \frac{|p|^j\ e^{\nu|p|}}{M_0(1+|p|^2)}\ \|F\|_\nu\ \|G\|_\nu \qquad (24)$$

PROOF. Let $p = |p|\ e^{i\theta}$. Then,

$$|(p^j\ F) * G(p,t)| \ =\ \left|\int_0^p \tilde{s}^j\ F(\tilde{s})\ G(p-\tilde{s})\ d\tilde{s}\right|$$

$$\leq \int_0^{|p|} ds\ s^j\ |F(s\ e^{i\theta})|\ |G(p - se^{i\theta})| \qquad (25)$$

Using the definition of $\|\cdot\|_\nu$, the above is bounded by

$$\frac{|p|^j}{M_0^2} e^{\nu|p|} \|F\|_\nu\ \|G\|_\nu \int_0^{|p|} \frac{ds}{(1+s^2)[1+(|p|-s)^2]}$$

$$\leq\ \frac{|p|^j\ e^{\nu|p|}}{M_0(1+|p|^2)}\ \|F\|_\nu\ \|G\|_\nu$$

The last inequality follows from the fact that

$$\int_0^{|p|} \frac{1}{(1+s^2)[1+(|p|-s)^2]}\ =\ 2\ \frac{\ln\ (|p|^2+1)\ +\ |p|\ \tan^{-1}\ |p|}{|p|\ (|p|^2\ +\ 4)}$$

and the definition of $M_0$ in (19). □

**Corollary 12** *With the convolution $*$, $\mathcal{A}_\phi$ is a Banach Algebra and furthermore,*

$$\|F * G\|_\nu\ \leq\ \|F\|_\nu\ \|G\|_\nu \qquad (26)$$

PROOF. This follows by applying Lemma 11 for $j=0$ and using the definition of $\|\cdot\|_\nu$.
□

**Lemma 13** *For $\nu > 2\rho_0 + 1$,*

$$\left|\int_0^t (p^j F) * B_{j,k} * F^{*k}\ e^{-p^3(t-\tau)}\ d\tau\right|$$

$$\leq\ \frac{C(\phi)}{M_0(1+|p|^2)}\ \||B_{j,k}| * |F|\ \|_\nu\ \|F\|_\nu^k\ e^{\nu|p|}\ T^{(3-j)/3} \qquad (27)$$

*where the constant $C$ is independent of $T$, but depends on $\phi$.*

PROOF. For $k \geq 1$, from Lemma 11, with $G = (B_{j,k} * F) * F^{*(k-1)}$, we obtain

$$|(p^j F) * B_{j,k} * F^{*k}| \leq \frac{|p|^j e^{\nu |p|}}{M_0(1+|p|^2)} \|B_{j,k} * F\|_\nu \, \|F\|_\nu^k$$

$$\leq \frac{|p|^j e^{\nu |p|}}{M_0(1+|p|^2)} \||B_{j,k}| * |F|\|_\nu \, \|F\|_\nu^k \quad (28)$$

For $k = 0$, we note that

$$|(p^j F) * B_{j,0}| \leq |p|^j \, \||F| * |B_{j,0}|\| \leq \frac{|p|^j e^{\nu |p|}}{M_0(1+|p|^2)} \, \||B_{j,0}| * |F|\|_\nu$$

Therefore, (28) holds for $k = 0$ as well. Thus, for any $k \geq 0$,

$$\left| \int_0^t (p^j F) * B_{j,k} * F^{*k} \, e^{-p^3(t-\tau)} \, d\tau \right|$$

$$\leq \frac{|p|^j e^{\nu |p|}}{M_0(1+|p|^2)} \, \||B_{j,k}| * |F|\|_\nu \, \|F\|_\nu^k \int_0^t e^{-|p|^3 \cos(3\theta)(t-\tau)} \, d\tau$$

where $p = |p| e^{i\theta}$. On integrating with respect to $\tau$, we obtain

$$|p|^j \int_0^t e^{-|p|^3 \cos(3\theta)(t-\tau)} \, d\tau \leq \frac{T^{(3-j)/3}}{\cos^{j/3} 3\phi} \sup_\gamma \frac{1 - e^{-\gamma^3}}{\gamma^{3-j}}$$

□

**Definition 14** For $H \in \mathcal{H}$, $F$ and $h$ in $\mathcal{A}_\phi$, define $h_0 = 0$ and for $k \geq 1$,

$$h_k \equiv H * [(F+h)^{*k} - F^{*k}]. \quad (29)$$

**Lemma 15** For $\nu > \rho + 1$,

$$\|h_k\|_\nu \leq k \Big( \|F\|_\nu + \|h\|_\nu \Big)^{k-1} \||H| * |h|\|_\nu \quad (30)$$

PROOF. We prove this by induction. $k = 0$ follows trivially since $h_0 = 0$. The case of $k = 1$ is obvious from (29). Assume the formula (30) holds for all $k \leq l$. Then

$$\|h_{l+1}\|_\nu = \|H * (F+h) * (F+h)^{*l} - H * F * F^{*l}\|_\nu = \|H * h * (F+h)^{*l} + F * h_l\|_\nu$$

On using (30) for $k = l$, we get

$$\leq \||H| * |h|\|_\nu (\|F\|_\nu + \|h\|_\nu)^l + l\|F\|_\nu \Big( \|F\|_\nu + \|h\|_\nu \Big)^{l-1} \||H| * |h|\|_\nu$$

$$\leq (l+1) \Big( \|F\|_\nu + \|h\|_\nu \Big)^l \||H| * |h|\|_\nu$$

Thus the formula (30) holds for $k = l + 1$. □

**Lemma 16** *For $F$ and $h$ in $\mathcal{A}_\phi$, and $\nu > 2\rho_0 + 1$, and $k \geq 1$,*

$$\left| (p^j[F+h]) * B_{j,k} * \left(F+h\right)^{*k} - (p^j F) * B_{j,k} * F^{*k} \right|$$
$$\leq \frac{|p|^j \, e^{\nu|p|}}{M_0(1+|p|^2)} \left( \|F\|_\nu + \|h\|_\nu \right)^{k-1} \Big\{ k\|F\|_\nu \||B_{j,k}| * |h|\|_\nu +$$
$$\|B_{j,k} * (F+h)\|_\nu \|h\|_\nu \Big\} \quad (31)$$

PROOF. It is clear that

$$\left| (p^j[F+h]) * B_{j,k} * \left(F+h\right)^{*k} - (p^j F) * B_{j,k} * F^{*k} \right|$$
$$\leq \left| (p^j h) * B_{j,k} * \left(F+h\right)^{*k} \right| + |(p^j F) * h_k| \quad (32)$$

where $H$ is now replaced by $B_{j,k}$ in the definition of $h_k$ in (29). Applying Lemma 11 and Corollary 12 to the first term, we obtain for $k \geq 1$

$$\left| (p^j h) * B_{j,k} * \left(F+h\right)^{*k} \right|$$
$$\leq \frac{|p|^j e^{\nu|p|}}{M_0(1+|p|^2)} \|h\|_\nu \|B_{j,k} * (F+h)\|_\nu \left( \|F\|_\nu + \|h\|_\nu \right)^{k-1}$$

On the other hand, applying Lemma 11 and Lemma 15, with $H = B_{j,k}$ and $\rho$ replaced by $2\rho_0$, we obtain

$$|(p^j F) * h_k| \leq \frac{|p|^j e^{\nu|p|}}{M_0(1+|p|^2)} k\|F\|_\nu \|B_{j,k} * h\|_\nu [\|F\|_\nu + \|h\|_\nu]^{k-1}$$

Combining the previous two equations, and using it in (32), we obtain the proof of Lemma 16 by noting that $\|B_{j,k} * h\|_\nu \leq \||B_{j,k}| * |h|\|_\nu$. □

**Lemma 17** *For $\nu > 2\rho_0 + 1$,*

$$\left| \int_0^t \left\{ (p^j[F+h]) * B_{j,0} - (p^j F) * B_{j,0} \right\} e^{-p^3(t-\tau)} \, d\tau \right|$$
$$\leq \frac{C(\phi) T^{(3-j)/3} \, e^{\nu|p|}}{M_0(1+|p|^2)} \||B_{j,0}| * |h|\|_\nu \quad (33)$$

PROOF. We note that

$$\left| (p^j[F+h]) * B_{j,0} - (p^j F) * B_{j,0} \right| \leq |(p^j h) * B_{j,0}|$$
$$\leq \frac{|p|^j e^{\nu|p|}}{M_0(1+|p|^2)} \||B_{j,0}| * |h|\|_\nu$$

Further, as before in the proof of Lemma 13

$$\int_0^t |p|^j \left| e^{-p^3(t-\tau)} \right| d\tau \leq C(\phi) T^{(3-j)/3}$$

Combining the two equations above, the Lemma follows.
□

**Lemma 18** *For $\nu > 2\rho_0 + 1$, and $k \geq 1$,*

$$\left| \int_0^t \left\{ (p^j[F+h]) * B_{j,k} * (F+h)^{*k} - (p^j F) * B_{j,k} * F^{*k} \right\} e^{-p^3(t-\tau)} d\tau \right|$$

$$\leq \frac{C(\phi) T^{(3-j)/3} e^{\nu|p|}}{M_0(1+|p|^2)} (\|F\|_\nu + \|h\|_\nu)^{k-1} \left\{ k\|F\|_\nu \||B_{j,k}| * |h|\|_\nu \right.$$

$$\left. + \|B_{j,k} * (F+h)\|_\nu \|h\|_\nu \right\} \quad (34)$$

PROOF. The proof is similar to that of Lemma 13, except Lemma 16 is used instead of Lemma 11.
□

**Lemma 19** *For $F \in \mathcal{A}_\phi$, and $\nu > 4\rho_0 + 1$ large enough so that $(\nu/2)^{-\beta}\|F\|_\nu < 1$, $\mathcal{N}F$ satisfies the following bounds*

$$\|\mathcal{N}F\|_\nu \leq A_b(T) \sum_{j=0}^3 C_j(\phi) T^{(3-j)/3} (\nu/2)^{-\alpha_j} \frac{\|F\|_\nu}{1 - (\nu/2)^{-\beta}\|F\|_\nu} + \|F_0\|_\nu \quad (35)$$

*Further, for $h \in \mathcal{A}_\phi$ such that $(\nu/2)^{-\beta}(\|F\|_\nu + \|h\|_\nu) < 1$*

$$\|\mathcal{N}(F+h) - \mathcal{N}F\|_\nu$$

$$\leq A_b(T) \sum_{j=0}^3 C_j(\phi)(\nu/2)^{-\alpha_j} T^{(3-j)/3} \frac{\|h\|_\nu}{[1 - (\nu/2)^{-\beta}(\|F\|_\nu + \|h\|_\nu)]^2} \quad (36)$$

PROOF. On inspection of (16) and using Lemma 13, it follows that

$$\|\mathcal{N}F\|_\nu \leq \sum_{j=0}^3 C_j(\phi) T^{(3-j)/3} \sum_{k=0}^\infty \||B_{j,k}| * |F|\|_\nu \|F\|_\nu^k + \|F_0\|_\nu \quad (37)$$

Now using Corollary 10 and noting the dependence of $C_1$ on $j$ through $\alpha_j$, (35) follows. As far as (36), from inspection of (16) and using Lemmas 17 and 18, we get

$$\|\mathcal{N}(F+h) - \mathcal{N}F\|_\nu \leq \sum_{j=0}^3 C_j(\phi) T^{(3-j)/3} \left( \| |B_{j,0}| * |h| \|_\nu \right.$$

$$\left. + \sum_{k=1}^\infty (\|F\|_\nu + \|h\|_\nu)^{k-1} \left\{ k\|F\|_\nu \||B_{j,k}| * |h|\|_\nu + \|B_{j,k} * (F+h)\|_\nu \|h\|_\nu \right\} \right) \quad (38)$$

On using Corollary 10 and noting the dependence of $C_1$ on $j$ through $\alpha_j$, (36) follows. □

**Comment 8:** Lemma 19 is the key to showing the existence and uniqueness of solution in $\mathcal{A}_\phi$ to (16), since it provides the conditions for the nonlinear operator $\mathcal{N}$ to map a ball into itself as well the necessary contractivity condition.

**Lemma 20** *If there exists some $b > 1$ so that*

$$(\nu/2)^{-\beta} b \|F_0\|_\nu < 1$$

$$\text{and} \quad A_b(T) \sum_{j=0}^{3} \frac{C_j(\phi)(\nu/2)^{-\alpha_j} T^{(3-j)/3}}{1 - (\nu/2)^{-\beta} b \|F_0\|_\nu} < 1 - \frac{1}{b} \quad (39)$$

*then the nonlinear mapping $\mathcal{N}$ maps a ball of radius $b\|F_0\|_\nu$ back into itself. Further, if*

$$A_b(T) \sum_{j=0}^{3} \frac{C_j(\phi)(\nu/2)^{-\alpha_j} T^{(3-j)/3}}{[1 - (\nu/2)^{-\beta} b \|F_0\|_\nu]^2} < 1 , \quad (40)$$

*then $\mathcal{N}$ is a contraction there.*

PROOF. This is a simple application of Lemma 19, if we simply note that $\|F\|_\nu^k < b^k \|F_0\|_\nu^k$. □

**Lemma 21** *For any given $T > 0$ and $\phi$ in the interval $(0, \pi/6)$, for all sufficiently large $\nu$, there exists a unique $F \in \mathcal{A}_\phi$ that satisfies the integral equation (16).*

PROOF.
We choose $b = 2$. It is clear from the bounds on $\|F_0\|_\nu$ in Lemma 7 that for given $T$, since $\alpha_r \geq 1$, $b(\nu/2)^{-\beta} \|F_0\|_\nu < 1$ for all sufficiently large $\nu$. Further, it is clear on inspection that both conditions (39) and (40) are satisfied for all sufficiently large $\nu$. The lemma now follows from the fixed point theorem.
□

## 4.1 Behavior of solution $F_s$ near $p = 0$

**Proposition 22** *For some $K_1 > 0$ and small $p$ we have $|F_s| < K_1 |p|^{\alpha_r - 1}$ and thus $|f_s| < K_2 |y|^{-\alpha_r}$ for some $K_2 > 0$ as $|y| \to \infty$ in $\mathcal{D}_{\rho,\phi}$.*

PROOF. Convergence in $\|\cdot\|_\nu$ implies uniform convergence on compact subsets of $\mathcal{K}$ and we can interchange summation and integration in (16). With $F_s$ the unique solution of (16) we let

$$G_j = \sum_{k=0}^{\infty} (-1)^j B_{j,k} * F_s^{*k}$$

and define the linear operator $\mathcal{G}$ by

$$\mathcal{G}Q = \int_0^t e^{-p^3(t-\tau)} \sum_{j=0}^{3}(p^j Q) * G_j \, d\tau$$

Clearly $F_s$ also satisfies the linear equation

$$F_s = \mathcal{G}F_s + F_0 \quad \text{or} \quad F_s = (1-\mathcal{G})^{-1}F_0$$

For $a > 0$ small enough define $\overline{\mathcal{S}}_a = \overline{\mathcal{S}} \cap \{p : |p| \leq a\}$. Since $F_s$ is continuous in $\overline{\mathcal{S}}$ we have $\lim_{a \downarrow 0} \|\mathcal{G}\| = 0$, where the norm is taken over $C(\overline{\mathcal{S}}_a)$.

By (2), (7), (17) and Lemma 3, we see that $\|F_0\|_\infty \leq K_3 |a|^{\alpha_r - 1}$ in $\overline{\mathcal{S}}_a$ for some $K_3 > 0$ independent of $a$. Then, as $a \downarrow 0$, we have

$$\max_{\overline{\mathcal{S}}_a} |F(p,t)| = \|F\| \leq (1 - \|\mathcal{G}\|)^{-1} \max_{\overline{\mathcal{S}}_a} \|F_0\| \leq 2K_3 |a|^{\alpha_r - 1}$$

and thus for small $p$ we have $|F(p,t)| \leq 2K_3 |p|^{\alpha_r - 1}$ and the proposition follows.

$\square$

**Proof of Theorem 2.** Lemma 3 shows that if $f$ is any solution of (1) satisfying Condition 1, then $\mathcal{L}^{-1}\{f\} \in \mathcal{A}_{\phi-\delta}$ for $0 < \delta < \phi$ for $\nu$ sufficiently large. For large $y$, the series (3) converges uniformly and thus $F = \mathcal{L}^{-1}\{f\}$ satisfies (16), which by Lemma 21 has a unique solution in $\mathcal{A}_\phi$ for any $\phi \in (0, \pi/6)$. Conversely, if $F_s \in \mathcal{A}_{\tilde{\phi}}$ is the solution of (16) for $\nu > \nu_1$, then from comment 6, $f_s = \mathcal{L}F_s$ is analytic in $\mathcal{D}_{\phi,\rho}$ for $0 < \phi < \tilde{\phi} < \pi/6$, for sufficiently large $\rho$, where in addition from Proposition 22, $f_s = O(y^{-\alpha_r})$. This implies that the series in (1) converges uniformly and by properties of Laplace transform, $f_s$ solves (1) and satisfies condition (1).

**Comment 9:** Theorem 2 can be applied directly to each of the examples 1-3 in the following sections to give existence and uniqueness of solution in $\mathcal{D}_{\phi,\rho}$ for any given time $T$, provided $\rho$ is large enough. However, this general theorem does not provide the specific dependence of $\rho$ on $T$. In the following sections, we not only show that Theorem 2 can be applied to the examples given, but use the specific information on $b_{j,k}(y,t)$ and $r(y,t)$ to obtain dependence of $\rho$ on $T$. This requires additional case-specific lemmas and theorems.

## 5 Example 1

This example comes in the context of solving the leading order inner-equation for a complex singularity of the conformal mapping function corresponding to an evolving Hele-Shaw flow ([12], Equations 5.5-5.9) as well as a 2-D dendrite in the small Peclet number limit, when surface energy anisotropy is small ([17], Equation A44 after a transformation). Consider the PDE

$$H_t = H^3 H_{xxx} \tag{41}$$

with the initial condition

$$H(x,0) = x^\gamma \tag{42}$$

where $0 < \gamma < 1$ (Note: this $\gamma$ is related to the $\beta$ defined in [12], through $\gamma = \beta/2$), and the far field matching condition

$$H(x,t) = x^\gamma + O(x^{4\gamma-3}) \tag{43}$$

as $|x| \to \infty$ for $\arg x \in (-\frac{2\pi}{3(1-\gamma)}, \frac{2\pi}{3(1-\gamma)})$. Here $\gamma$ is real and in the interval $(0,1)$.

The transformations:

$$y = \frac{x^{1-\gamma}}{1-\gamma} \quad ; \quad H = x^\gamma(1+y^{-1}f(y,t)) \tag{44}$$

bring the equation to the form (1), with

$$r(y,t) = -\gamma(\gamma-1)^{-2}(\gamma-2)y^{-2} \tag{45}$$

and the only nonzero coefficients are $\{b_{i,j} : i,j = 0,1,2,3\}$, are presented below as a matrix.

$$\begin{bmatrix} \frac{22\gamma - 11\gamma^2 - 6}{y^3(\gamma-1)^2} & 9\frac{6\gamma - 3\gamma^2 - 2}{y^4(\gamma-1)^2} & \frac{50\gamma - 25\gamma^2 - 18}{y^5(\gamma-1)^2} & 2\frac{(1-2\gamma)(2\gamma-3)}{y^6(\gamma-1)^2} \\ \frac{7\gamma^2 - 14\gamma + 6}{(\gamma-1)^2 y^2} & 3\frac{7\gamma^2 - 14\gamma + 6}{y^3(\gamma-1)^2} & 3\frac{7\gamma^2 - 14\gamma + 6}{y^4(\gamma-1)^2} & \frac{7\gamma^2 - 14\gamma + 6}{y^5(\gamma-1)^2} \\ -3y^{-1} & -9y^{-2} & -9y^{-3} & -3y^{-4} \\ 0 & 3y^{-1} & 3y^{-2} & y^{-3} \end{bmatrix}$$

The initial condition (42) translates as

$$f_I(y) = 0 \tag{46}$$

Also, the condition (43) implies that as $y \to \infty$ for $(y,t) \in \mathcal{D}_{\phi,\rho}$ (as defined earlier),

$$f(y,t) = O(y^{-2}) \tag{47}$$

From the expressions for $b_{j,k}$ above it is possible to calculate $B_{j,k}$ explicitly. For our purpose, it is enough to note that aside from $B_{3,0}$, which is identically zero, we can write

$$|B_{j,k}(p,t)| < C|p|^{2-j+k}$$

for a constant $C$ independent of $j$ and $k$, as well as $T$ and therefore from Lemma 7, we conclude that for $F \in \mathcal{A}_\phi$, for $(j,k) \neq (3,0)$,

$$\||B_{j,k}| * |F|\|_\nu < C \nu^{-3+j-k} \|F\|_\nu \tag{48}$$

We also note that from Lemma7 that

$$\|F_0(p,t)\|_\nu < A_r \nu^{-1} T \tag{49}$$

where $A_r$ here is independent of $T$. Since only a finite number of $B_{j,k}$ are nonzero, it is better to use properties (37) and (38) directly to obtain conditions for the fixed point theorem to apply:

$$\frac{1}{b} + C \sum_{j=0}^{3}{}' \sum_{k=0}^{3} b^k \ A_r^k \ T^k \nu^{-2k+j-3} T^{(3-j)/3} \ < \ 1 \tag{50}$$

$$\sum_{j=0}^{3}{}' \sum_{k=0}^{3} b^k \ (k+1) \ A_r^k \ T^k \nu^{-2k+j-3} T^{(3-j)/3} \ < \ 1 \tag{51}$$

Here the primes in the summation symbol in (50) and (51) mean that the term $j = 3$, $k = 0$ is missing. Each of these conditions (50) and (51) are satisfied for $\nu^{-1} T^{1/3}$ sufficiently small for any choice of $b > 1$. The condition that $T/\nu^3$ is less than some number translates into $T/\rho^3$ being sufficiently small, i.e. we are restricted to a region of space in the $x$-plane where $T/x^{3(1-\gamma)}$ is sufficiently small. It was noted earlier [Kadanoff, private communication, 1991 and independently by Howison [16] for a special case] that there was a similarity solution to (41) of the form:

$$H(x,t) = t^{\gamma/(3(1-\gamma))} \ h\left(x/t^{\frac{1}{3(1-\gamma)}}\right) \tag{52}$$

The resulting ordinary differential equation for $h(\eta)$ was solved numerically [12] and a first few singularities (in order of distance from the origin) of $H$ were thus determined. It was surmised that these solutions form an infinite set that straddles the boundary of the sector $|\arg x| < \frac{2\pi}{3(1-\gamma)}$, over which one can specify the asymptotic condition $h(\eta) \sim \eta^\gamma$ that one needs to satisfy initial and far-field conditions (42-43). Later, these conclusions were confirmed rigorously by Fokas and Tanveer [14], who transformed the equation into Painlevé $P_{II}$ and used isomonodromic approaches for integrable systems to confirm the earlier behavior seen numerically. These results are also amenable to exponential asymptotic methods and formal trans-series association with actual functions, that have recently been worked out [3]. The latter method is more general than the isomonodromic methods since they apply equally well for integrable and nonintegrable equations. The application of our existence and uniqueness results mean that the only solution to the initial value problem for the PDE are those with similarity structure given by (52).

From what has been discussed so far and proved in this paper, an interesting aspect of the complex plane initial value problem is that that the initial condition (42) is not recovered as $t \to 0^+$, except in the sector $|\arg x| < \frac{2\pi}{3(1-\gamma)}$. This follows from the equivalence of small $t$ with large $x$ in the similarity solution (52).

## 6 Example 2

The second example also comes from Hele-Shaw flow ([12], equations (6.10)-(6.12)) as well as dendritic crystal growth for weak undercooling and for weak anisotropic surface energy [17]. In the asymptotic limit of surface tension tending to zero, it was determined that in the neighborhood of an initial zero, the local governing equation is

$$H_t \ + \ H_x \ = \ H^3 \ H_{xxx} - \frac{1}{2} H^3 \tag{53}$$

The initial condition is:

$$H(x, 0) = x^{-1/2} \tag{54}$$

The far-field matching condition with the "outer" asymptotic solution is

$$H(x, t) = x^{-1/2} + O(x^{-5}) \tag{55}$$

as $|x| \to \infty$ with arg $x \in (-4\pi/9, 4\pi/9)$. It is to be noted that in this case, unlike case I, there are no similarity solutions which satisfy both the initial and asymptotic boundary conditions. We introduce the transformation:

$$x = t + \left(\frac{3y}{2}\right)^{2/3}, \quad H(x, t) = x^{-1/2} + x^{-3/2} y^{-1} f(y, t) \tag{56}$$

Note that if $x$ is large enough, $y \sim \frac{2}{3} x^{3/2}$. The initial condition (54) implies that

$$f_I(y) = 0 \tag{57}$$

and the asymptotic far-field condition (55) implies:

$$f(y, t) = O(y^{-4/3}) \tag{58}$$

as $y \to \infty$ in some $\mathcal{D}_{\phi,\rho}$. Under the change of variables, the PDE is of the form (1) with

$$r(y, t) = -\frac{15y}{8x^{7/2}} \tag{59}$$

and each $b_{j,k}$ containing one or more terms of the form $x^{-\beta} y^{-\delta}$, with $\beta > 0$ and $\frac{2}{3}\beta + \delta > 0$. The exact expressions for $b_{j,k}$ are given in the appendix.

**Comment 10:** Since the conditions for Theorem 2 hold, we may simply apply it and obtain a unique analytic solution over a sector $\mathcal{D}_{\phi,\rho}$ for a fixed $\phi$ satisfying conditions $0 < \phi < \pi/6$. The theorem, when applied to (53), would imply that for any $T$, there exists unique analytic solution $H(x, t)$ in the sector arg $x \in (-\pi/3 - 3\phi/2, -\pi/3 + 3\phi/2)$ provided $|x|$ is large enough. Since this is true for any $\phi$ in the interval $(0, \pi/6)$, the theorem establishes the existence of a unique analytic solution assumed before [12]. However, the restriction on how large $x$ has to be depends on $T$ and because of the generality, Theorem 2 does not give a precise dependence on $T$. Finding this constraint is the objective of the rest of this section.

**Lemma 23** *Let $g(y, t) = x^{-\beta} y^{-\delta}$, where $x$ is given by (56), $\beta > 0$ and $\frac{2}{3}\beta + \delta > 0$, then for any $p \in \mathcal{S}_\phi$,*

$$|G(p, t)| \leq C \, p^{\frac{2}{3}\beta + \delta - 1} \tag{60}$$

*where $C$ is independent of $\nu$ and $T$, but can depend on $\beta$ and $\delta$. Also, if in addition, $\frac{2}{3}\beta + \delta > 1$, then for $\nu > 1$,*

$$\|G\|_\nu \leq C \, \nu^{-2\beta/3 - \delta + 1} \tag{61}$$

PROOF. Given the relation between $x$ and $y$, we note that we may write

$$g(y, t) = y^{-2\beta/3-\delta} h(ty^{-2/3})$$

where

$$h(s) = \left(\frac{3}{2}\right)^{-2\beta/3} \left(1 + \frac{2^{2/3}}{3^{2/3}} s\right)^{-\beta}$$

It is to be noted that for arg $s$ bounded away from $\pm\pi$, $h(s)$ is uniformly bounded. It is also clear that

$$G(p, t) = \mathcal{L}^{-1} g[p, t] = \frac{1}{2\pi i} p^{2\beta/3+\delta-1} \left(\int_C e^s h(tp^{2/3} s^{-2/3}) \, ds\right)$$

where the contour $C$ is similar to that shown in Fig. 1, except in the $s$-plane. The intersection point of $C$ on the real $s$-axis will be chosen to be 1. It is clear that for $p \in \mathcal{S}_\phi$ and $s$ on the contour $C$, $\arg(tp^{2/3} s^{-2/3}) \in (-5\pi/9, 5\pi/9)$. Therefore $h$ is uniformly bounded and therefore

$$|G(p, t)| < C |p|^{\frac{2}{3}\beta+\delta-1} \int_0^\infty e^{-r/2} dr$$

Hence the lemma follows. ☐

**Comment 11:** Note that the preceding Lemma gives a much sharper result than applying the more general Lemma 3, because we made specific use of the form of the function $g(y,t)$.

**Corollary 24** $\|R\|_\nu < C \nu^{-1/3}$.

PROOF. This follows simply from the expression for $r$ and and application of Lemma 23. ☐

**Corollary 25** *For some $C$ independent of $T$ and $p$,*
(i) $|B_{0,0}| < C |p|^2$, $|B_{0,1}| < C |p|^{4/3}$, $|B_{0,2}| < C |p|^3$ and $|B_{0,3}| < C |p|^7$
(ii) $|B_{1,0}| < C|p|$, $|B_{1,1}| < C|p|^{8/3}$, $|B_{1,2}| < C|p|^{13/3}$ and $|B_{1,3}| < C|p|^6$.
(iii) $|B_{2,0}| < C$, $|B_{2,1}| < C |p|^{5/3}$, $|B_{2,2}| < C |p|^{10/3}$ and $|B_{2,3}| < C |p|^5$.
(iv) $|B_{3,0}| < C T |p|^{-1/3}$, $|B_{3,1}| < C |p|^{2/3}$, $|B_{3,2}| < C |p|^{7/3}$ and $|B_{3,3}| < C |p|^4$.

PROOF. These inequalities follow immediately from the expressions of the coefficients $b_{j,k}$ in the appendix and Lemma 23. ☐

**Corollary 26** *For $C$ independent of $T$, $p$, $j$ and $k$ and $\nu > 1$,*

$$\||B_{j,k}| * |F|\|_\nu < C \nu^{2j/9-2/3-k/3} \|F\|_\nu \text{ for } (j,k) \neq (3,0)$$

$$\||B_{3,0}| * |F|\|_\nu < C T \nu^{-2/3} \|F\|_\nu$$

PROOF. The estimates follow immediately on examination of the upper-bounds on $B_{3,k}$ and using Lemma 9 with $\rho = 0$. ☐

**Comment 12:** Application of Lemma 9 to the results Corollary 25 leads to stronger bounds; however the bounds indicated in Corollary 26 suffice for our purposes.

**Lemma 27** *For any $0 < \phi < \pi/6$, the integral equation (16), with $B_{j,k}$ and $R$ as determined in this section, has a unique solution $F_s \in \mathcal{A}_\phi$ provided $T \nu^{-2/3} < \epsilon$, where $\epsilon$ (depending only on $\phi$) is small. Further, $F_s = O(p^{1/3})$ as $p \to 0$ in $\mathcal{S}_\phi$*

PROOF. Since $F_0(p,t)$ is given by (17), and $F_I(y) = 0$, it follows that $\|F_0\|_\nu < CT\nu^{-1/3}$. Further, from the estimates of Corollary 26, it follows from (37) that the condition for mapping a ball of radius $b\|F_0\|_\nu$ back into itself is now given by

$$\frac{1}{b} + C \sum_{j,k=0}^{3}{}' \left(T\nu^{-2/3}\right)^{(3-j)/3} \left(bT\nu^{-2/3}\right)^k + CT\nu^{-2/3} < 1$$

where $\sum'$ denotes is the summation without $(j,k) = (3,0)$ term. From (38), the contractivity requirement becomes

$$C \sum_{j,k=0}^{3}{}' (k+1)\, b^k \left(T\nu^{-2/3}\right)^{(3-j)/3} \left(bT\nu^{-2/3}\right)^k + CT\nu^{-2/3} < 1$$

It is clear that both conditions are satisfied for some $b > 1$ (say $b = 2$) if $T \nu^{-2/3}$ is chosen sufficiently small. This ensures existence of a solution in the Banach space $\mathcal{A}_\phi$. Since $C$ in the above equations depend on $\phi$, the upper bound of $T \nu^{-2/3}$ for which the solution is guaranteed to exist depends on $\phi$. Further, since $r(y,t) = O(y^{-4/3})$, by applying proposition 22, it follows that $F_s(p,t) = O(p^{1/3})$ as $p \to 0$ in $\mathcal{S}_\phi$. □

**Theorem 28** *For any $0 < \phi < \pi/6$, there exists a unique solution $H(x,t)$ satisfying (53)-(55) that is analytic in $x$ in the domain*

$$\left\{(x,t) : |x| > \tilde{\rho},\ \arg\, x \in \left(-\frac{4}{9}\pi + \frac{2}{3}\phi, \frac{4}{9}\pi - \frac{2}{3}\phi\right),\ 0 \le t \le T\right\}$$

*provided $T\tilde{\rho}^{-1} < \epsilon$, with $\epsilon$ (depending only on $\phi$) small enough. Thus, for any $t$, when $\arg\, x \in (-4\pi/9, 4\pi/9)$, there exists a unique analytic solution satisfying (53-55) when $t/|x|$ is small enough.*

PROOF. By applying the equivalence between the solution to the integral equation (16), $F_s \in \mathcal{A}_\phi$ for $0 < \phi < \pi/6$, to the analytic solution $f_s$ of the partial differential equation (1) satisfying condition 1 (as shown in the proof of theorem 2), it follows that in this particular example, a solution $f_s$ exists for $(y,t) \in \mathcal{D}_{\phi,\rho}$, provided $T \rho^{-2/3} < \epsilon$ and that $f_s = O(y^{-4/3})$ as $y \to \infty$ in $\mathcal{D}_{\phi,\rho}$. Theorem 28 follows merely from noting the change of variable $(y,t,f)$ to $(x,t,H)$, once we choose $\tilde{\rho} = \rho^{2/3}$ and use the relation $y \sim \frac{2}{3}x^{3/2}$ for large $x$. □

# 7  Example 3: Strongly anisotropic inner equation

For strong anisotropic surface energy, the analytically continued conformal mapping function that maps the upper half-plane to the exterior of a one-sided two dimensional dendritic interface for small Peclet number satisfies, upon transformation, the following leading order inner-equation near a singularity of particular type (See [17], equation (A16) after some elementary transformations):

$$H_t = H^{1/3} H_{xxx} \tag{62}$$

with initial and boundary conditions

$$H(x, 0) = x^{-9\delta} \tag{63}$$

$$H(x, t) = x^{-9\delta} + O(x^{-12\delta-3}) \text{ as } |x| \to \infty \tag{64}$$

for arg $x \in \left(-\frac{2\pi}{3(1+\delta)}, \frac{2\pi}{3(1+\delta)}\right)$ where $\delta > 0$. If we introduce the transformation

$$y = \frac{x^{\delta+1}}{1+\delta} \; ; \; H(x,t) = x^{-9\delta}(1 + f(y,t)/y) \tag{65}$$

then we obtain an equation of the form (1) with

$$r(y, t) = \frac{9\delta(9\delta + 1)(9\delta + 2)}{(\delta + 1)^3 y^2} \tag{66}$$

$$b_0(f, y, t) = \frac{9\delta(9\delta + 1)(9\delta + 2)\{(1 + f/y)^{4/3} - 1\}y}{(\delta + 1)^3 y^3 f}$$

$$+ (54\delta^2 + 277\delta + 32)\frac{(1 + f/y)^{1/3}}{(\delta + 1)^2 y^3} \tag{67}$$

$$b_1(f, y, t) = \left(-\frac{217\delta + 26}{(1 + \delta)^2} - \frac{48\delta}{\delta + 1} - 6\right)(1 + f/y)^{1/3} y^{-2} \tag{68}$$

$$b_2(f, y, t) = \frac{3(9\delta + 1)(1 + f/y)^{1/3}}{y(\delta + 1)} \tag{69}$$

$$b_3(f, y, t) = 1 - (1 + f/y)^{1/3} \tag{70}$$

Using the series expansions in $f/y$ it follows that

$$b_{0,k} = \frac{9\delta(9\delta + 1)(9\delta + 2)\binom{4/3}{k+1}}{(\delta + 1)^3 y^{4+k}} + \frac{(54\delta^2 + 277\delta + 32)\binom{1/3}{k+1}}{(\delta + 1)^2 y^{3+k}}$$

$$b_{1,k} = \left\{-\frac{217\delta + 26}{(1 + \delta)^2} - \frac{48\delta}{\delta + 1} - 6\right\}\binom{1/3}{k} y^{-k-2}$$

$$b_{2,k} = \frac{3(9\delta + 1)\binom{1/3}{k}}{y^{k+1}(\delta + 1)}$$

$$b_{3,0} = 0$$

$$b_{3,k} = -\binom{1/3}{k} y^{-k} \text{ for } k \geq 1$$

It is clear that

$$
\begin{aligned}
B_{0,k} &= \frac{9\delta(9\delta+1)(9\delta+2)\binom{4/3}{k+1}p^{k+3}}{(\delta+1)^3(k+3)!} + \frac{(54\delta^2+277\delta+32)\binom{1/3}{k}p^{k+2}}{(\delta+1)^2(k+2)!} \\
B_{1,k} &= \left\{-\frac{217\delta+26}{(1+\delta)^2} - \frac{48\delta}{\delta+1} - 6\right\}\binom{1/3}{k}\frac{p^{k+1}}{(k+1)!} \\
B_{2,k} &= \frac{3(9\delta+1)\binom{1/3}{k}p^k}{(\delta+1)k!} \\
B_{3,0} &= 0 \\
B_{3,k} &= -\binom{1/3}{k}\frac{p^{k-1}}{(k-1)!} \quad \text{for } k \geq 1
\end{aligned}
$$

Thus, the following estimates hold:

$$\||B_{j,k}|*|F|\|_\nu \;<\; C\nu^{-k+j-3}\,\|F\|_\nu$$

where $C$ is a constant that can be made independent of $j$, $k$, $T$. Therefore, using the above relation, from (37), the condition for mapping a ball of radius $b\,\|F_0\|_\nu$ back into itself becomes:

$$C\sum_{j=0}^{3}{}'\sum_{k=0}^{\infty}(\nu^{-3}T)^{(3-j)/3}\left(\nu^{-1}b\|F_0\|_\nu\right)^k + \frac{1}{b} \;<\; 1$$

where the $\sum'$ indicates that $j=3$, $k=0$ term is missing from the summation. Applying the estimates of this section to (38), the contraction condition is:

$$C\sum_{j=0}^{3}{}'\sum_{k=0}^{\infty}(\nu^{-3}T)^{(3-j)/3}(k+1)\left(\nu^{-1}b\|F_0\|_\nu\right)^k \;<\; 1$$

Since $\|F_0\|_\nu < KT\nu^{-1}$, a sufficient condition for use of contraction mapping theorem is that

$$bKT\nu^{-2} \;<\; 1$$

i.e. that $T\nu^{-2}$ is small enough. Note that in that case $T\nu^{-3}$ is automatically small when $\nu$ is sufficiently large. The restriction $T\nu^{-2}$ small means that the differential equation (62), with conditions (63- 64) has a unique analytic solution for any $x \in \left(-\frac{2\pi}{3(1+\delta)}, \frac{2\pi}{3(1+\delta)}\right)$ in a region where when $tx^{-2\delta-2}$ is small enough. However, (62) admits a similarity solution

$$H(x,t) \;=\; t^{-\frac{3\delta}{1+\delta}}\, q\left(x/t^{\frac{1}{3(1+\delta)}}\right)$$

and $q(\eta)$ solves an ordinary differential equation and the asymptotic boundary condition

$$q(\eta) \;\sim\; \eta^{-9\delta} \quad \text{for } \eta \to \infty,\; \arg\eta \in \left(-\frac{2\pi}{3(1+\delta)}, \frac{2\pi}{3(1+\delta)}\right)$$

Uniqueness means that the similarity solution is the only solution. However, the restriction $tx^{-2\delta-2} \ll 1$ is sub-optimal; from the similarity structure of the solution, one expects analyticity in a sector for $tx^{-3\delta-3} \ll 1$, i.e. the integral equation (16) should have a unique solution for $T\nu^{-3}$ sufficiently small. The condition for convergence of the infinite series involving the norm of $F$ is stronger than the original condition for the convergence of the infinite series in $f$ that appear in the expansion of $(1+f/y)^{1/3}$. This is the reason for the sub-optimal estimates in this example involving infinite series. Nonetheless, the uniqueness results even for a restricted range, prove that the similarity solution determined earlier is the only solution to the problem.

# 8  Conclusion

We have proved existence and uniqueness of solution to a class of third-order nonlinear partial differential equation in a sector of the complex spatial variable $y$ for sufficiently large $|y|$. Our technique is akin to Borel summation, that has demonstrated its effectiveness in the analysis of general classes of nonlinear ODEs [3]. The class of PDEs for which existence and uniqueness has been proved contains three examples that arise in the context of Hele-Shaw fluid flow and dendritic crystal growth. The uniqueness results shows that the similarity solutions assumed earlier for examples 1 and 3 are the only ones that satisfy the given initial and far-field matching conditions. Accordingly, the singularities of the PDE solutions are those that correspond to singularities of the similarity solutions.

Expressions for $b_{j,k}$ for Example 2

$$b_{0,0} = -\frac{35}{6}\frac{1}{x^{3/2}y^2} - \frac{75}{4}x^{-9/2} - \frac{45}{8}\frac{(12)^{1/3}}{x^{7/2}y^{2/3}} - \frac{15}{4}\frac{(18)^{1/3}}{x^{5/2}y^{4/3}}$$

$$b_{0,1} = -\frac{35}{2}\frac{1}{x^{5/2}y^3} - \frac{45}{x^{11/2}y} - \frac{3}{2x^2 y} - \frac{45}{4}\frac{(18)^{1/3}}{x^{7/2}y^{7/3}} - \frac{135}{8}\frac{(12)^{1/3}}{x^{9/2}y^{5/3}}$$

$$b_{0,2} = -\frac{165}{4}\frac{1}{x^{13/2}y^2} - \frac{135}{8}\frac{(12)^{1/3}}{x^{11/2}y^{8/3}} - \frac{1}{2x^3 y^2} - \frac{35}{2}\frac{1}{x^{7/2}y^4} - \frac{45}{4}\frac{(18)^{1/3}}{x^{9/2}y^{10/3}}$$

$$b_{0,3} = -\frac{45}{8}\frac{(12)^{1/3}}{x^{13/2}y^{11/3}} - \frac{35}{6}\frac{1}{x^{9/2}y^5} - \frac{105}{8}\frac{1}{x^{15/2}y^3} - \frac{15}{4}\frac{(18)^{1/3}}{x^{11/2}y^{13/3}}$$

$$b_{1,0} = \frac{15}{4}\frac{(18)^{1/3}}{x^{5/2}\sqrt[3]{y}} + \frac{45}{8}\frac{\sqrt[3]{y}(12)^{1/3}}{x^{7/2}} + \frac{35}{6}\frac{1}{x^{3/2}y}$$

$$b_{1,1} = \frac{35}{2}\frac{1}{x^{5/2}y^2} + \frac{45}{4}\frac{(18)^{1/3}}{x^{7/2}y^{4/3}} + \frac{135}{8}\frac{(12)^{1/3}}{x^{9/2}y^{2/3}}$$

$$b_{1,2} = \frac{35}{2}\frac{1}{x^{7/2}y^3} + \frac{135}{8}\frac{(12)^{1/3}}{x^{11/2}y^{5/3}} + \frac{45}{4}\frac{(18)^{1/3}}{x^{9/2}y^{7/3}}$$

$$b_{1,3} = \frac{35}{6}\frac{1}{x^{9/2}y^4} + \frac{15}{4}\frac{(18)^{1/3}}{x^{11/2}y^{10/3}} + \frac{45}{8}\frac{(12)^{1/3}}{x^{13/2}y^{8/3}}$$

$$b_{2,0} = -3\,x^{-3/2} - \frac{9(18)^{1/3}y^{2/3}}{4x^{5/2}}$$

$$b_{2,1} = -\frac{9}{x^{5/2}y} - \frac{27}{4}\frac{(18)^{1/3}}{x^{7/2}\sqrt[3]{y}}$$

$$b_{2,2} = -\frac{9}{x^{7/2}y^2} - \frac{27(18)^{1/3}}{4x^{9/2}y^{4/3}}$$

$$b_{2,3} = -\frac{3}{x^{9/2}y^3} - \frac{9(18)^{1/3}}{4x^{11/2}y^{7/3}}$$

$$b_{3,0} = -1 + \frac{3y}{2x^{3/2}}$$

$$b_{3,1} = \frac{9}{2}x^{-5/2}$$

$$b_{3,2} = \frac{9}{2x^{7/2}y}$$

$$b_{3,3} = \frac{3}{2x^{9/2}y^2}$$

**Acknowledgments** OC was supported in part by the National Science Foundation (NSF-DMS-9704968) as well as by the Ohio State University Math Research Institute, while ST acknowledges partial support from the National Science Foundation (NSF-DMS-9803358) and NASA (NAG 3-1947).